\documentclass[a4paper,notitlepage, 8pt,reqno]{article}

  \usepackage{mathtext}
       \usepackage{graphicx}
 \usepackage[cp1251]{inputenc}
  \usepackage[T2A]{fontenc}
 \usepackage[english]{babel}
 \usepackage[all,arc,poly,2cell,curve,arrow,tips]{xy}

  \usepackage{enumerate, eucal, amsthm,amsmath, amssymb}

  \allowdisplaybreaks[4]

 \theoremstyle{definition}
 \newtheorem{defn}{Definition}

 \theoremstyle{plain}
 \newtheorem{thm}{Theorem}

 \newtheorem{prop}{Proposion}

\newtheorem{cor}{Corrolary}

 \newtheorem{lem}{Lemma}

 \theoremstyle{remark}

 \renewcommand{\abstractname}{}

  \newcounter{ab}
\title{
Noncommutative pfaffians and representations. Applications to
classification of states of five-dimensional quasi-spin. \footnote
 {The work is supported by grants  Nsh-8508.2010.1,  MK-1378.2009.1}}
\author{Dmitry Artamonov\footnote{artamonov.dmitri@gmail.com},  Valentina Goloubeva\footnote{goloubeva@yahoo.com}}

\begin{document}

\maketitle

\renewcommand{\abstractname}{}

\begin{abstract} Noncommutative pfaffians associated with an orthogonal  algebra are some
 special elements of the universal enveloping algebra. Is is proved that  in the case $N=2n+1$ some
 of these pfaffians, denoted as $PfF_{\widehat{-n}}$ and  $PfF_{\widehat{n}}$, act on the space of $\mathfrak{o}_{2n-1}$-highest vectors
of a $\mathfrak{o}_{2n+1}$-representation.  There exist the
Mickelsson-Zhelobenko algebra of raising operators
$Z(\mathfrak{o}_{2n+1},\mathfrak{o}_{2n-1})$ which naturally acts on
this space. We find explicitly an element of the
Mickelsson-Zhelobenko algebra , which acts on this space in the same
way as the pfaffian $PfF_{\widehat{n}}$. As a by product we find
explicit formulae for the action of the pfaffian $PfF_{\widehat{n}}$
  in the Gelfand-Tsetlin-Molev base.
 The action
of pfaffians in the tensor realization of
 representation is considered in the appendix.

\end{abstract}

\maketitle{}

\section{Introduction}
In the paper we study the noncommutative pfaffians, which are some
special elements of the universal enveloping  algebra
$U(\mathfrak{o}_N)$.  The main subject is an investigation of an
action of noncommutative pfaffians in representations in the case of
odd $N$.

First  it  is noted that all pfaffians for  map weight vectors to
weight vectors and the weight changes by a simple rule.

Then it is noted that in the case $N=2n+1$ some pfaffian (these
pfaffians are denoted as $PfF_{\widehat{-n}}$ and
$PfF_{\widehat{n}}$, the explanation of these notations see below in
the introduction) commute with the subalgebra
$\mathfrak{o}_{2n-1}\subset \mathfrak{o}_{2n+1}$. Thus these
pfaffians act on the space of $\mathfrak{o}_{2n-1}$-highest vectors
of a $\mathfrak{o}_{2n+1}$ representation.

Thus the pfaffians $PfF_{\widehat{-n}}$, $PfF_{\widehat{n}}$ act as
 raising operators in the problem of the construction of a bases of
a $\mathfrak{o}_{2n+1}$-representation of Gelfand-Tsetlin type. The
problem is to construct a bases in a
$\mathfrak{o}_{2n+1}$-representation and for different
  $n$ the constructions  must be coherent.
  This means that  the bases of a $\mathfrak{o}_{2n+1}$-representation
  must
  be a union of bases in $\mathfrak{o}_{2n-1}$-representation into
   which the $\mathfrak{o}_{2n+1}$-representation splits.



There  papers where such a base is constructed in the simplest
nontrivial case $\mathfrak{o}_3\subset\mathfrak{o}_5$
(\cite{Goli},\cite{Rore1},\cite{Rore2},\cite{Rore3}).

In  general case, using ideas of Gelfand and Tsetlin, such a base
was constructed by Molev \cite{Mol1},  \cite{Mol3}, \cite{Mol4}, see
also \cite{Molev1} and  \cite{Molev}. This base is called in the
present paper the Gelfand-Tsetlin-Molev base. Mention that a basis
of the same type in the case $\mathfrak{sp}_{2n}$ was firstly
constructed by Zhelobenko \cite{zb}.

In the construction of the the Gelfand-Tsetlin-Molev  base the key
role is played by the Mickelsson-Zhelobenko algebra
$Z(\mathfrak{o}_N,\mathfrak{o}_{N-2})$.
  (see it's definitions in sec. \ref{mj}), which acts on the space of $\mathfrak{o}_{N-2}$-highest
  vectors of a $\mathfrak{o}_N$-rerpresentation. The elements of this
  algebra are called raising operators.

The first main result of the paper is the following. We find
explicitly an element of
$Z(\mathfrak{o}_{2n+1},\mathfrak{o}_{2n-1})$, which acts on the
space of $\mathfrak{o}_{2n-1}$-highest vectors of a
$\mathfrak{o}_{2n+1}$-representation in the same way as the pfaffian
$PfF_{\widehat{n}}$.

This element is proportional to  $\check{z}_{n0}$ which is one of
the canonical generators of the Mickelsson-Zhelobenko algebra. A
coefficient of proportionality belongs to
$U(\mathfrak{h}_{\mathfrak{o}_{2n-1}})$.

One can interpret the first main result as follows. Each operator
acting on the space of $\mathfrak{o}_{2n-1}$-highest vectors of a
$\mathfrak{o}_{2n+1}$-representation can be continued to an
operator, acting on the whole
  $\mathfrak{o}_{2n+1}$-representation. This continuation commutes with the action of
    $\mathfrak{o}_{2n-1}$.
We find this continuation for the element  $C\check{z}_{n0}$ of the
Mickelsson-Zhelobenko,
  $C\in U(\mathfrak{h}_{\mathfrak{o}_{2n-1}})$. For other raising operator such continuation
  are not known  in an explicit way as far as we know.

The second main result are explicit formulae for the action of
$PfF_{\widehat{n}}$ in the Gelfand-Tsetlin-Molev base (Theorem
\ref{maint}).

Let us give the main definitions and describe the structure of the
paper.

Let  $\Phi=(\Phi_{ij})$, $i,j=1,...,2n$  be a skew-symmetric
 $2n\times 2n$-matrix, whose  matrix entries belong to
  a noncommutative ring.

\begin{defn}The noncommutative pfaffian
 of $\Phi$ is defined by the formulae

$$Pf \Phi=\frac{1}{n!2^n}\sum_{\sigma\in
S_{2n}}(-1)^{\sigma}\Phi_{\sigma(1)\sigma(2)}...\Phi_{\sigma(2n-1)\sigma(2n)},$$

Here $\sigma$ is a permutation of the set  $\{1,...,2n\}$.
\end{defn}

In the paper a split realization of $\mathfrak{o}_N$ is used. The
algebra $\mathfrak{o}_{N}$ is generated by matrices
$F_{ij}=E_{ij}-E_{-j-i}$, here $E_{ij}$ are matrix units. When $N$
is odd the indices $i,j$ belong to the set
$\{-n,...,-1,0,1,...,n\}$, where $n=\frac{N-1}{2}$. When $N$ is even
the indices $i,j$ belong to the set $\{-n,...,-1,1,...,n\}$, where
$n=\frac{N}{2}$. Shortly in both cases this set of indices is
denoted in the paper as $\{-n,...,n\}$.

In the paper the following noncommutative pfaffians are  considered.

\begin{defn} Let $F$ be the matrix $F=(F_{ij})$.
For every subset $I\subset\{-n,...,n\}$ which consists of an even
number $k$ of elements define a submatrix $F_I$   by the formulae
$F_{I}=(F_{ij})_{-i,j\in I}$. Put

$$PfF_I=Pf(F_{-ij})_{-i,j\in I}.$$
\end{defn}

In \cite{Molev}  the author  in terms of these pfaffians defines
some special elements of $U(\mathfrak{o}_N)$ called the Capelli
elements. These elements are
 $C_k=\sum_{I\subset \{-n,...,n\},|I|=k}PfF_IPfF_{-I}$,
$k=2,4,...,[\frac{N}{2}]$.  It is proved that the elements $C_k$
belong to the center of $U(\mathfrak{o}_N)$.

Mention that the relation beween commutative pfaffians were
intensively studied in \cite{Dress},
\cite{Ish1},\cite{Ish2},\cite{Ish3}. Also some relation between
noncommutative pfaffians $PfF_I$  were derived in
\cite{I},\cite{AG}.

The structure if the paper is the following.

In   sec.  \ref{wei} the action of pfaffians on  weight vectors is
investigated. It is proved that a weight vector is mapped to a
weight vector, a rule according to which the weight changes is
established (Proposition \ref{wes}).

In sec. \ref{or} some facts about the split realization of the
orthogonal algebra are given. In sec. \ref{com}, it is proved that
in the case $N=2n+1$ the pfaffians
$$PfF_{\widehat{-n}}:=PfF_{\{-n+1,...,n\}} \text{ and } PfF_{\widehat{n}}:=PfF_{\{-n,...,n-1\}}$$ commute with elements of
the subalgebra $\mathfrak{o}_{2n-1}=<F_{ij}>$, $-n+1\leq i,j \leq
n-1$ (Corrolary \ref{cor3}).

In particular in the case $\mathfrak{o}_{2n+1}$ the pfaffian
$PfF_{\widehat{-n}}$ diminishes the $n$-th component of the weight
by one and $PfF_{\widehat{n}}$ raises the $n$-th component of the
weight by one (Corollary \ref{wes1}).

In sec \ref{formu} some formulaes involving pfaffians are proved.

In  sec. \ref{repr} using two the facts mentioned above an important
observation is done.  The pfaffians  $PfF_{\widehat{-n}}$  and
$PfF_{\widehat{n}}$ act on the space of
$\mathfrak{o}_{2n-1}$-highest vectors of a
$\mathfrak{o}_{2n+1}$-representation  with a fixed
$\mathfrak{o}_{2n-1}$-weight.

Then we find an element of the the Mickelsson-Zhelobenlo algebra
that acts on the space of $\mathfrak{o}_{2n-1}$-highest vectors in
the same way as the pfaffian $PfF_{\widehat{n}}$.

There exists a projection from $U(\mathfrak{o}_{N})$  to
$Z(\mathfrak{o}_N,\mathfrak{o}_{N-2})$. The element of the algebra
that acts in the same way as the pfaffian  is an image of
$PfF_{\widehat{n}}$ under this projection. After some calculations
in sec. \ref{mjp1}, \ref{mjp2} we do it (Theorem \ref{tm}).

Using this result in    sec. \ref{dpb}  we find   explicit formulae
for the action of the pfaffian $PfF_{\widehat{n}}$ on the base in
the Gelfand-Tsetlin-Molev base (Theorem \ref{maint}).

Application of these calculations are presented in Sec. \ref{quasi}.
We formulate a problem of construction of an additional quantum
number for classifications of quasi-spin states and establish its
relation to the Gelfans-Tsetlin-Molev bases for
$\mathfrak{o}_5$-representations.

In subsection \ref{number} Theorem \ref{teork} is proved,  in which
we construct an additional quantum number using pfaffians.

 In     appendix   an action of a pfaffian in a tensor  representation is investigated.
  The base  vectors are encoded by orthogonal Young tableaus \cite{proc}.
   An action of a pfaffian on  tensor products of vectors of a standard representation  is found
(Propositions  \ref{stan},\ref{er},\ref{Pfk},\ref{Pft}). Then
theorem \ref{Pfvt} is proved which gives  an information about the
action on base vectors given by Young tableaus. The image of a base
vector in this theorem is expressed  as a linear combination of not
necessarily orthogonal Young tableaus. Thus Theorem \ref{Pfvt}
 does not give explicit formulaes of the action of a pfaffian in the
bases formed by orthogonal Young tableaus.

\section{The split realization of the orthogonal algebra}
\label{or}

The orthogonal algebra $\mathfrak{o}_N$ is a tangent space at the
unit to the group of linear transformations that preserve a
nondegenerate quadratic form. Let
 $G$ be a matrix of the form. A matrix $f$ belongs to the algebra
  $\mathfrak{o}_N$  if  $f^tG+Gf=0$.  We use the following indexation of
  rows  and colunmes of matrices $f$. When $N$
is odd the indices $i,j$ belong to the set
$\{-n,...,-1,0,1,...,n\}$, where $n=\frac{N-1}{2}$. When $N$ is even
the indices $i,j$ belong to the set $\{-n,...,-1,1,...,n\}$, where
$n=\frac{N}{2}$. Shortly this set of indices is denoted in the paper
as $\{-n,...,n\}$.

 In the paper the split realization
  of  the algebra  $\mathfrak{o}_N$  is used. It corresponds
   to the following choice of the form: $G=(\delta_i^{-j})$.

In this realization the algebra $\mathfrak{o}_N$ is generated by
elements $F_{ij}=E_{ij}-E_{-j-i}$. The only linear relations between
these elements are $F_{ij}=-F_{-j-i}$.

Commutation relations between these elements are

$$[F_{ij},F_{kl}]=\delta_{kj}F_{il}-\delta_{il}F_{kj}-\delta_{-ki}F_{-jl}
+\delta_{-lj}F_{k-i}.$$

One can prove that elements
$F_{-n-n},...,F_{-1-1}$ form a bases in the Cartan subalgebra, and  elements
$F_{ij},j<-i$ are root elements.

The precise  the correspondence is the following. Let $e_i$ be the
element  $F_{ii}^*$ in the dual space to the Cartan subalgebra. Put
$e_{-r}:=-e_{r}$ and $e_0=0$.  Then the element  $F_{ij}$
corresponds to the root $e_i-e_{j}$ ($\S 13$ in the chapter 8 in
\cite{Bu}).

\section{Noncommutative pfaffians.}

In this section some properties of noncommuttive pfaffians are
obtained. Firstly the action  on weight vectors of
 a representation is investigated. Then commutators of
 pfaffians with elements $F_{ij}$ are calculated.  Finally some summation formulaes involving
 pfaffians are proved.

\subsection{Action of a pfaffian on a weight vector.}
\label{wei}
 Remind that $e_i$ denotes the standard base vectors $F_{ii}^*$ of  dual space to the Cartan subalgebra.

\begin{prop}
 \label{wes}
Let $V$ be a representation of $\mathfrak{o}_N$.
  Under the action of the pfaffian $PfF_I$  a weight vector  with the weight $\mu$ is mapped to a weight vector with the weight $\mu-\sum_{i\in I}e_i$.
 \end{prop}
\proof
 If  $v$ is a weight vector in a representation  of $\mathfrak{o}_N$ with the weight $\mu$, $g_{\alpha}$ is a root element in $\mathfrak{o}_N$ corresponding to the root $\alpha$,
 then $g_{\alpha}v_{\mu}$ is a weight vector of to the weight  $\alpha+\mu$.

Consider the vector $PfF_Iv$.
  By definition one has
  $$PfF_I=\frac{1}{\frac{k}{2}!2^{\frac{k}{2}}}\sum_{\sigma\in
S^k}(-1)^{\sigma}F_{-\sigma(i_1)\sigma(i_2)}...F_{-\sigma(i_{k-1})\sigma(i_k)}.$$

To prove the proposition it suffices to show that every summand
changes the weight by substracting of the same expression
$-\sum_{i\in I}e_i$. Using the correspondence between roots and
elements $F_{ij}$ from the sec. \ref{or} one gets the following.
When one acts by
$F_{-\sigma(i_1)\sigma(i_2)}...F_{-\sigma(i_{k-1})\sigma(i_k)}$ on
$v$ then to the weight the vector
$$e_{-\sigma(i_1)}-e_{\sigma(i_2)}-...+e_{-\sigma(i_{k-1})}-e_{\sigma(i_k)}=-\sum_{i\in
I}e_i$$ is added. This proves the proposition.

 \endproof

Consider the most interesting case
$\mathfrak{o}_{N}=\mathfrak{o}_{2n+1}$ and $|I|=2n$.

\begin{cor}
\label{wes1} Let $\mathfrak{o}_{N}=\mathfrak{o}_{2n+1}$.

The action of $PfF_{\widehat{-n}}$ adds the vector
 $-\sum_{i\in I}e_i=-e_{n}$ to the weight.

The action of $PfF_{\widehat{n}}$ adds the vector
 $-\sum_{i\in I}e_i=-e_{-n}=e_{n}$ to the weight.
\end{cor}

\subsection{Commutators of pfaffians and $F_{ij}$.}

\label{com}

\begin{lem}
\label{l2}
  Let $I=\{i_1,...,i_k\}$, where $k$ is even. Then the commutator $[PfF_I,F_{j_1-j_2}]$
   is calculated according to the following rule.

\begin{enumerate}

\item If $j_1,j_2\notin I$, then  $[PfF_I,F_{j_1-j_2}]=0$.

\item If $j_1\in I,j_2\notin I$, then $[PfF_I,F_{j_1-j_2}]=PfF_{I_{j_1\rightarrow -j_2}}$.

\item If $j_1\notin I,j_2\in I$, then  $[PfF_I,F_{j_1-j_2}]=-PfF_{I_{j_2\rightarrow -j_1}}$.

\item If $j_1\in I,j_2\in I$, then
  $[PfF_I,F_{j_1-j_2}]=PfF_{I_{j_1\rightarrow -j_2}}-PfF_{I_{j_2\rightarrow -j_1}}$
\end{enumerate}

\end{lem}

\proof

One can identify $E_{ij}$ with $e_{i}\otimes e_{-j}$. Then $F_{ij}$
is identified with $e_{i}\wedge e_{-j}$. Remind that
 $$Pf
F_I=\frac{1}{2^{\frac{k}{2}}\frac{k}{2}!}\sum_{ \sigma\in
S_k}(-1)^{\sigma}F_{-\sigma(i_1)\sigma(i_2)}...F_{-\sigma(i_{k-1})\sigma(i_{k})}.$$

Thus  $PfF_I$ with indexing set $I=\{i_1,...,i_k\}$ is  identified
with the polyvector $e_{-i_1}\wedge...\wedge e_{-i_k}$.

This identification is compatible with the action of
$\mathfrak{o}_N$. Thus

\begin{center}
$[PfF_I,F_{j_1-j_2}]=-[F_{j_1-j_2},PfF_I]=F_{j_1-j_2}e_{i_1}\wedge
e_{i_2}...\wedge e_{i_k}+e_{i_1}\wedge F_{j_1-j_2}e_{i_2}... \wedge
e_{i_k}+e_{i_1}\wedge e_{i_2}... \wedge
F_{j_1-j_2}e_{i_k}.$\end{center}

One has $F_{j_1-j_2}e_{-i_1}=e_{j_2}$ if $j_2=i_1$ and
$F_{j_1-j_2}e_{-i_1}=-e_{j_2}$ if $j_1=i_1$.

Suggest that $\{j_1,j_2\}\cap I=\emptyset$.   Then
$[PfF_I,F_{j_1-j_2}]=0$.

Suggest that $j_1\in I,j_2\notin I$. Then  $j_1=i_t$ and the only
nonzero summand is that containing $F_{j_1j_2}e_{-i_t}$.  Thus we
have $[PfF_I,F_{j_1-j_2}]=-PF_{I_{j_1\mapsto -j_2}}$. Here
$I_{j_1\mapsto -j_2}$ is obtained from $I$ by replacing the index
$j_1$ to $j_2$.

The case $j_1\notin I,j_2\in I$ is considered in the same manner.

Suggest that $\{j_1,j_2\}\subset I$. That is $j_1=e_{i_t}$,
$j_2=e_{i_s}$. Then in the sum there are two nonzero summands one
contains  $F_{j_1j_2}e_{i_t}$, the other contains
$F_{j_1j_2}e_{i_s}$.  Each of them is a wedge product with a new
indexing set. So one gets that
$[PfF_I,F_{j_1j_2}]=PfF_{I_{j_1\rightarrow
-j_2}}-PfF_{I_{j_2\rightarrow -j_1}}$.

\endproof

Write the formulae from Lemma \ref{l2} in the case $N=2n+1$ for
pfaffians $$PfF_{\widehat{i}}:=PfF_{\{-n,...,\widehat{i},...,n\}}.$$


\begin{cor}
\label{cor1} Let $N=2n+1$, then
\begin{enumerate}
\item $[PfF_{\widehat{-i}},F_{i-j}]=(-1)^{i+j}PfF_{\widehat{j}}.$
\item $[PfF_{\widehat{-j}},F_{i-j}]=-(-1)^{i+j}PfF_{\widehat{i}}.$
\item If $k\neq i,j$
$[PfF_{\widehat{k}},F_{i-j}]=0.$

\end{enumerate}
\end{cor}

\proof Take formulaes for $[PfF_{I},F_{i-j}]=0$, written in Lemma
\ref{l2}.

If  $i\notin I$ or $j\notin I$,  then the considered commutator is
zero.

If $i,j\in I$, than the considered commutator equals to
$PfF_{I_{i\rightarrow -j}}-PfF_{I_{j\rightarrow -i}}$. Consider this
case in details.

 Note that $I$ contains all indices except one.

 If $-i,-j\in I$ (that is $PfF_{I}\neq
PfF_{\widehat{-i}},PfF_{\widehat{-j}}$), then both  summands
$PfF_{I_{i\rightarrow -j}},PfF_{I_{j\rightarrow -i}}$ vanish.

If $-i\notin I$, then $PfF_{I_{i\rightarrow -j}}=0$,
$PfF_{I_{j\rightarrow -i}}=PfF_{\{-n,...,-i_{\text{at the place
j}},...,\widehat{-i},...,n\}}=(-1)^{i+j+1}PfF_{\widehat{j}}$. That
is in this case
$[PfF_{\widehat{-i}},F_{i-j}]=(-1)^{i+j}PfF_{\widehat{j}}$.

If $-j\notin I$, then $PfF_{I_{j\rightarrow -i}}=0$,
$PfF_{I_{i\rightarrow -j}}=PfF_{\{-n,...,-j_{\text{at the place
i}},...,\widehat{-j},...,n\}}=(-1)^{i+j+1}PfF_{\widehat{i}}$. That
is in this case
$[PfF_{\widehat{-j}},F_{i-j}]=(-1)^{i+j+1}PfF_{\widehat{i}}$.

\endproof

\begin{cor}
\label{cor3} In the case $\mathfrak{o}_{2n+1}$  the pfaffians
$PfF_{\widehat{n}}$, $PfF_{\widehat{-n}}$ commute with elements
$F_{ij}$, $-n<i,j<n$, that span the subalgebra $\mathfrak{o}_{2n-1}$

\end{cor}

\subsection{Some formulas involving pfaffians.}
\label{formu}

In this subsection some  summation formulaes  are
proved.

\begin{lem}
\label{minors1}

$PfF_I=\frac{(\frac{p}{2})!(\frac{q}{2})!}{(\frac{k}{2})!}\sum_{I=I'\sqcup
I'',|I'|=p,|I''|=q}(-1)^{(I'I'')}PfF_{I'}PfF_{I''}$.

Here $(-1)^{(I'I'')}$ is a sign of a permutation of the set
$I=\{i_1,...,i_k\}$ that places first the subset $I'\subset I$ and
then the subset $I''\subset I$.

The numbers $p$, $q$ are even fixed numbers,  they satisfy
$p+q=k=|I|$.
\end{lem}
\proof

By definition one has

$$PfF_I=\frac{1}{2^{\frac{k}{2}}(\frac{k}{2})!}\sum_{\sigma\in
S_k}(-1)^{\sigma}F_{-\sigma(i_1)\sigma(i_2)}...F_{-\sigma(i_{k-1}),\sigma(i_k)}.$$

The summand $
(-1)^{\sigma}F_{-\sigma(i_1)\sigma(i_2)}...F_{-\sigma(i_{k-1}),\sigma(i_k)}$
can be written as

$$(-1)^{(I'I'')}(-1)^{\sigma'}F_{-\sigma'(i'_1)\sigma'(i'_2)}...F_{-\sigma'(i'_{p-1}),\sigma'(i'_p)}
(-1)^{\sigma''}F_{-\sigma''(i''_1)\sigma''(i''_2)}...F_{-\sigma''(i''_{q-1}),\sigma''(i''_q)}.$$

Here $I'=\{i'_1,...,i'_p\}$ is the set of indices
$\{\sigma(i_1),...,\sigma(i_p)\}$ placed in a natural order,
$I''=\{i''_1,...,i''_q\}$ is set   of indices
$\{\sigma(i_{p+1}),...,\sigma(i_k)\}$ placed in a natural order,
$\sigma'$ is a permutation $\{\sigma(i_1),...,\sigma(i_p)\}$  of the
set $I'$ and  $\sigma''$ is a permutation of the set $I''$ defined
in a similar way. Note that
$(-1)^{(I'I'')}(-1)^{\sigma'}(-1)^{\sigma''}=(-1)^{\sigma}$.

The mapping $\sigma\mapsto I',I'',\sigma',\sigma''$ is bijective.

Thus the pfaffian can be written as

$\frac{(\frac{p}{2})!(\frac{q}{2})!}{(\frac{k}{2})!}\sum_{I=I'\sqcup
I'',|I'|=p,|I''|=q}(-1)^{(I'I'')}\frac{1}{2^{\frac{k}{2}}(\frac{p}{2})!(\frac{q}{2})!}\sum_{\sigma'}(-1)^{\sigma'}(-1)^{\sigma''}
F_{-\sigma'(i'_1)\sigma'(i'_2)}...\\...F_{-\sigma'(i'_{p-1}),\sigma'(i'_p)}F_{-\sigma''(i''_1)\sigma'(i''_2)}...F_{-\sigma''(i''_{q-1}),\sigma''(i''_q)}=\\
=\frac{(\frac{p}{2})!(\frac{q}{2})!}{(\frac{k}{2})!}\sum_{I=I'\sqcup
I'',|I'|=p,|I''|=q}(-1)^{(I'I'')}PfF_{I'}PfF_{I''}$
\endproof

\begin{cor}
\label{corl2} If $|I|=k$, then

$PfF_I=\frac{1}{\frac{k}{2}+1}\sum_{I=I'\sqcup
I''}(-1)^{(I'I'')}\frac{(\frac{|I'|}{2})!(\frac{|I''|}{2})!}{(\frac{k}{2})!}PfF_{I'}PfF_{I''}$.

\end{cor}

\begin{lem}
\label{minorn}

 Let $-n\in I$. Then
$PfF_I=\\=\sum_{i\in I\setminus\{-n\}}\sum_{I\setminus
\{-n,i\}=I'\sqcup
I''}\frac{(\frac{|I'|}{2})!(\frac{|I''|}{2})!}{(\frac{k}{2})!}(-1)^{(I'-niI'')}PfF_{I'}F_{ni}PfF_{I''}$.

Here $(-1)^{(I'-niI'')}$ is a sign of the permutation $(I',-n,i,I'')$
of the set $I$.
\end{lem}
\proof

By definition one has

$$PfF_I=\frac{1}{2^{\frac{k}{2}}(\frac{k}{2})!}\sum_{\sigma\in
S_k}(-1)^{\sigma}F_{-\sigma(i_1)\sigma(i_2)}...F_{-\sigma(i_{k-1}),\sigma(i_k)}.$$

Since $F_{ij}=-F_{-j-i}$ the summation can be taken only over such
permutation  such that $\sigma(i_{2t-1})<\sigma(i_{2t})$. But if the
summation is done in such a way the multiple
$\frac{1}{2^{\frac{k}{2}}}$ must be omitted.

Fix a such a permutation $\sigma$ and find a place such that
$(\sigma(i_{2t-1}),\sigma(i_{2t}))=(-n,i)$. The summand
$(-1)^{\sigma}F_{-\sigma(i_1)\sigma(i_2)}...F_{-\sigma(i_{k-1}),\sigma(i_k)}$
can be written as

$(-1)^{(I'inI'')}(-1)^{\sigma'}F_{-\sigma'(i'_1)\sigma'(i'_2)}...F_{-\sigma'(i'_{p-1}),\sigma'(i_p)}F_{-ni}(-1)^{\sigma''}F_{-\sigma''(i''_1)\sigma''(i''_2)}...
F_{-\sigma''(i''_{q-1}),\sigma''(i''_q)}$.

Here $I'=\{i'_1,...,i'_p\}$ is the set of indices
$\{\sigma(i_1),...,\sigma(i_{2t-2})\}$ placed in the natural order,
$I''=\{i''_1,...,i''_q\}$ is set   of indices
$\{\sigma(i_{2t+1}),...,\sigma(i_k)\}$ placed in the natural order,
$\sigma'$ is a permutation $\{\sigma(i_1),...,\sigma(i_{2t-2})\}$ of
the set $I'$ and $\sigma''$ is a permutation of the set $I''$
defined in a similar way. Note that
$(-1)^{(I'-niI'')}(-1)^{\sigma'}(-1)^{\sigma''}=(-1)^{\sigma}$. The
permutation $\sigma'$  satisfies the condition
$\sigma'(i'_{2t-1})<\sigma'(i'_{2t})$ as well as  the permutation
$\sigma''$.

The mapping $\sigma\mapsto I',I'',\sigma',\sigma''$ is bijective
(since $\sigma(i_{2t-1})<\sigma(i_{2t})$).

Thus the pfaffian can be written as

$\sum_{i\in I\setminus\{n\}}\sum_{I\setminus \{i,n\}=I'\sqcup I''}
\frac{(\frac{|I'|}{2})!(\frac{|I''|}{2})!}{(\frac{k}{2})!}\frac{1}{(\frac{|I'|}{2})!(\frac{|I''|}{2})!}(-1)^{(I'-niI'')}
(-1)^{\sigma'}F_{-\sigma'(i'_1)\sigma'(i'_2)}...\\...F_{-\sigma'(i'_{p-1}),\sigma'(i_p)}F_{ni}(-1)^{\sigma''}F_{-\sigma''(i''_1)\sigma''(i''_2)}...
F_{-\sigma''(i''_{q-1}),\sigma''(i''_q)}=\\=\sum_{i\in
I\setminus\{n\}}\sum_{I\setminus \{i,n\}=I'\sqcup
I''}\frac{(\frac{|I'|}{2})!(\frac{|I''|}{2})!}{(\frac{k}{2})!}(-1)^{(I'-niI'')}PfF_{I'}F_{ni}PfF_{I''}$

\endproof

\begin{lem}
\label{comult}

$\Delta PfF_I=\sum_{I'\sqcup I''=I}(-1)^{(I'I'')}PfF_{I'}\otimes
PfF_{I''}$

Here $(-1)^{(I'I'')}$ is a sign of a permutation of the set
$I=\{i_1,...,i_k\}$ that places first the subset $I'\subset I$ and
then places the subset $I''\subset I$.
\end{lem}
\proof

By definition one has

$$PfF_I=\frac{1}{2^{\frac{k}{2}}(\frac{k}{2})!}\sum_{\sigma\in
S_k}(-1)^{\sigma}F_{-\sigma(i_1)\sigma(i_2)}...F_{-\sigma(i_{k-1}),\sigma(i_k)}.$$

Apply the comuptilication, one gets

$$\Delta PfF_I=\frac{1}{2^{\frac{k}{2}}(\frac{k}{2})!}\sum_{\sigma\in
S_k}(-1)^{\sigma}(F_{-\sigma(i_1)\sigma(i_2)}\otimes 1+1\otimes
F_{-\sigma(i_1)\sigma(i_2)})...(F_{-\sigma(i_{k-1}),\sigma(i_k)}\otimes
1+1\otimes F_{-\sigma(i_{k-1}),\sigma(i_k)}) $$

The product $$(F_{-\sigma(i_1)\sigma(i_2)}\otimes 1+1\otimes
F_{-\sigma(i_1)\sigma(i_2)})...(F_{-\sigma(i_{k-1}),\sigma(i_k)}\otimes
1+1\otimes F_{-\sigma(i_{k-1}),\sigma(i_k)})$$ equals to
$$\sum_{I=J'\sqcup J''} F_{-\sigma(j'_1)\sigma(j'_2)}...
F_{-\sigma(j'_{p-1})\sigma(j'_p)}\otimes
F_{-\sigma(j''_1)\sigma(j''_2)}...
F_{-\sigma(j''_{q-1})\sigma(j'_q)}.$$ Here
$J'=\{j'_1,j'_2,...,j'_p\}$, $J''=\{j''_1,j''_2,...,j''_q\}$ are
subset of $I$, such that $I=J'\sqcup J''$ and also the following
condition are satisfied. If $\sigma(i_{2t-1})\in J'$ then
$\sigma(i_{2t})\in J'$, if $\sigma(i_{2t-1})\in J''$ then
$\sigma(i_{2t})\in J''$. In other word the partitions $I=J'\sqcup
J''$ must induce a division of $\frac{k}{2}$ pairs
$(\sigma(i_{2t-1}),\sigma(i_2t))$

The summand
$$(-1)^{\sigma}F_{-\sigma(j'_1)\sigma(j'_2)}...
F_{-\sigma(j'_{p-1})\sigma(j'_p)}\otimes
F_{-\sigma(j''_1)\sigma(j''_2)}...
F_{-\sigma(j''_{q-1})\sigma(j'_q)}$$ indexed by $\sigma\in Aut(I),
J',J''$ can be written as the following expression

$$(-1)^{(I'I'')}(-1)^{\sigma'}F_{-\sigma'(i'_1)\sigma'(i'_2)}...
F_{-\sigma'(i'_{p-1})\sigma'(i'_p)}\otimes
(-1)^{\sigma''}F_{-\sigma''(i''_1)\sigma''(i''_2)}...
F_{-\sigma''(i''_{p-1})\sigma''(i''_p)}$$

indexed by $I',I''$, $\sigma'\in Aut (I')$, $\sigma''\in Aut (I'')$.
Here $I'$ is the set $\{\sigma(j'_1),...,\sigma(j'_p)\}$ written in
the natural order, $I''$ is the set
$\{\sigma(j''_1),...,\sigma(j''_q)\}$ written in a natural order.
The permutation $\sigma'$ is the permutation
$\{\sigma(j'_1),...,\sigma(j'_p)\}$ of $J'$ and $\sigma''$ is a
permutation of $J''$ defined in a similar way.

But the mapping $\sigma, J',J''\mapsto I',I'', \sigma' , \sigma''$
is not injective.   To get the triple $\sigma, J',J''$ with  the
prescribed image $I',I'', \sigma' , \sigma''$ one must divide
$\frac{k}{2}$ pairs $\{(i_1,i_2),...,(i_{k-1},i_k)\}$ into two
subsets $J'$ and $J''$ with $\frac{|I'|}{2}$ and $\frac{|I''|}{2}$
elements respectively. Take a permutation $\sigma$, such that
$\sigma(J')=\sigma'(I')$ (as order sets), and
$\sigma(J'')=\sigma''(I'')$  (as order sets).  The only freedom is
the choice of two subsets $J'$ and $J''$. Thus the number of
elements in the preimage equals to the number divisions of
$\frac{k}{2}$ pairs into two subsets: one consists of
$\frac{|I'|}{2}$ pairs  and the other consists of $\frac{|I''|}{2}$
pairs. The number is
$\frac{(\frac{k}{2})!}{(\frac{|I'|}{2})!(\frac{|I''|}{2})!}$.

Thus $\Delta PfF_I$ can be
written as

$\frac{1}{2^{\frac{k}{2}}(\frac{k}{2})!}\sum_{I=I'\sqcup
I''}\frac{(\frac{k}{2})!}{(\frac{|I'|}{2})!(\frac{|I''|}{2})!}(-1)^{(I'I'')}(\sum_{\sigma'\in
Aut(I') }(-1)^{\sigma'}F_{-\sigma'(i'_1)\sigma'(i'_2)}...
F_{-\sigma'(i'_{p-1})\sigma'(i'_p)})\otimes\\ \otimes
(\sum_{\sigma''\in Aut(I'')
}(-1)^{\sigma''}F_{-\sigma''(i''_1)\sigma''(i''_2)}...
F_{-\sigma''(i''_{p-1})\sigma''(i''_p)})$.

This expression equals $\sum_{I'\sqcup
I''=I}(-1)^{(I'I'')}PfF_{I'}\otimes PfF_{I''}$.

\endproof

\section{Pfaffians are raising operators}
\label{repr}


According to the sections \ref{com} and \ref{wei}  the following holds

\begin{enumerate}
\item  The action of  $PfF_{\widehat{n}}$ and  $PfF_{\widehat{-n}}$ commutes with
the action of the subalgebra  $\mathfrak{o}_{2n-1}$, spanned by $F_{ij}$, $-n<i,j<n$.

\item  Under the action of  $PfF_{\widehat{-n}}$  a weight vector is mapped to a weight vector, the $n$-th component of the weight is dimished by $1$.
Under the  action of $PfF_{\widehat{n}}$  a weight vector is also
mapped to a weight vector, the $n$-th component of the weight is
raised by $1$.
\end{enumerate}

The following lemma is proved.

\begin{lem}
\label{hw} The pfaffians  $PfF_{\widehat{n}}$ and
$PfF_{\widehat{-n}}$ act on the space of  $\mathfrak{o}_{2n-1}$-highest vectors of a $\mathfrak{o}_{2n+1}$-representation. The
$\mathfrak{o}_{2n-1}$-weight under this action is conserved.

\end{lem}

\section{The Mickelson-Zhelobenko algebra and pfaffians.}\label{mj}

There exists the Mickelsson-Zhelobenko algebra which as
$PfF_{\widehat{n}}$ and $PfF_{\widehat{-n}}$ acts on the space of
  $\mathfrak{o}_{2n-1}$ highest vectors of a   $\mathfrak{o}_{2n+1}$ representation.
 In the present section  we  find an element of this  algebra which acts as $PfF_{\widehat{n}}$.

At first in the subsection \ref{mj1} we give the definition of the
 Mickelsson-Zhelobenko algebra. There exists a mapping from
$U(\mathfrak{o}_N)$ to the Mickelson-Zhelobenko algebra
$Z(\mathfrak{o}_{2n+1},\mathfrak{o}_{2n-1})$. The image of the
pfaffian $PfF_{\widehat{n}}$ is exactly an element of
$Z(\mathfrak{o}_{2n+1},\mathfrak{o}_{2n-1})$, which acts on the
space of  $\mathfrak{o}_{2n-1}$-highest vectors as the pfaffian.

 In
subsections  \ref{mjp1}, \ref{mjp2} the images of some special
pfaffain in $Z(\mathfrak{o}_{2n+1},\mathfrak{o}_{2n-1})$ are found.
Using this calculations at the end of the sunsection \ref{mjp2} the
image of  $PfF_{\widehat{n}}$ in the Mickelsson-Zhelobenko algebra
is found.

\subsection{The Mickelson-Zhelobenko algebra}\label{mj1}

The Gelfand-Tsetlin-Molev's approach to  a construction of a bases
of a $\mathfrak{o}_{2n+1}$-representation is based on restrictions
 $\mathfrak{o}_{2n+1}\downarrow \mathfrak{o}_{2n-1}$, in contrast to the classical Gelfand-Tsetlin's
  approach which is  based on restrictions $\mathfrak{o}_N\downarrow \mathfrak{o}_{N-1}$.
The subalgebra $\mathfrak{o}_{2n-1}\subset \mathfrak{o}_{2n+1}$ is
spanned by the elements $F_{ij}$,  $-n<i,j<n$. The Cartan subalgebra
$h_{\mathfrak{o}_{2n-1}}$ is a subalgebra in
$h_{\mathfrak{o}_{2n+1}}$ and root vectors in $\mathfrak{o}_{2n-1}$
are also root vectors in $\mathfrak{o}_{2n+1}$.

Remind a scheme of construction of a
$\mathfrak{o}_{2n+1}$-representation $V$. An irreducible
representation $V$ of the algebra $\mathfrak{o}_{2n+1}$ becomes
reducible as a representation of $\mathfrak{o}_{2n-1}$. According to
the scheme of Gelfand and Tsetlin in order  to construct a base one
firstly must know possible highest weights $\mu$  of irreducible
$\mathfrak{o}_{2n-1}$-representations into which splits $V$.
Secondly if a weight has multiplicity one must be able to construct
a bases in the  multiplicity space, that is in the space of
$\mathfrak{o}_{2n-1}$-highest vectors with a fixed
$\mathfrak{o}_{2n-1}$-weight $\mu$.

Introduce a notation for this space.

\begin{defn}
Let $V_{\mu}^{+}$ be a space of $\mathfrak{o}_{2n-1}$-highest
vectors with the $\mathfrak{o}_{2n-1}$-weight $\mu$ in a
$\mathfrak{o}_{2n+1}$-representation $V$.
\end{defn}

To construct a base in $V_{\mu}^{+}$   Molev used the
Mickelson-Zhelobenko algebra acting on the space
$\oplus_{\mu}V_{\mu}^+$.

Let us give a definition of this algebra, see also
 \cite{JL},\cite{Molev1},  and the chapter 9 in \cite{Molev}.

Let $\mathfrak{g}$ be a Lie algebra and let $\mathfrak{k}$ be it's
reductive subalgebra. The main example is
$\mathfrak{g}=\mathfrak{o}_{2n+1}$ and
$\mathfrak{k}=\mathfrak{o}_{2n-1}$. Let
$\mathfrak{k}=\mathfrak{k}^-+\mathfrak{h}+\mathfrak{k}^+$ be a
triangular decomposition. Let $R(\mathfrak{h})$ be a field of
fractions of the algebra $U(\mathfrak{h})$. Denote as
$$U'(\mathfrak{g})=U(\mathfrak{g})\otimes_{U(\mathfrak{h})}
R(\mathfrak{h}).$$

Let $$J'=U'(\mathfrak{g})\mathfrak{k}^+$$ be the left ideal in
$U'(\mathfrak{g})$, generated by $\mathfrak{k}^+$. Put
$$M(\mathfrak{g},\mathfrak{k})=U'(\mathfrak{g})/J'.$$

For every positive root $\alpha$ of the algebra $\mathfrak{k}$
define

$$p_{\alpha}=1+\sum_{k=1}^{\infty}e_{-\alpha}^ke_{\alpha}^k\frac{(-1)^k}{k!(h_{\alpha}+\rho(h_{\alpha})+1)...(h_{\alpha}+\rho(h_{\alpha})+k)},$$

here $e_{\alpha}$ is a root vector  $\mathfrak{k}$, corresponding to
$\alpha$, $h_{\alpha}$ is a corresponding Cartan element, $\rho$ is
a half-sum of positive roots of $\mathfrak{k}$.

An order is normal if the following holds. Let a root
  be
   a sum of two roots, then it lies between them. Chose a normal ordering $\alpha_1<...<\alpha_m$  of  positive roots
of
  $\mathfrak{k}$.

Put $$p=p_{\alpha_1}...p_{\alpha_m}.$$ This element is called the
extremal projector. It can be proved that nevertheless  $p$ is an
infinite series it's action on $M(\mathfrak{g},\mathfrak{k})$  by
left multiplication is well defined \cite{JL}.

The following equalities hold: $e_{\alpha}p=pe_{-\alpha}=0$, here $\alpha$  is a positive root of $\mathfrak{k}$.

Put $$Z(\mathfrak{g},\mathfrak{k})=pM(\mathfrak{g},\mathfrak{k}).$$
This is the Mickelson-Zhelobenko algebra. The multiplication in
$Z(\mathfrak{g},\mathfrak{k})$ is defined using the isomorphism
$Z(\mathfrak{g},\mathfrak{k})=Norm J'/J'$, where $NormJ'=\{u\in
U'(\mathfrak{g}): J'u\subset J'\}$. Thus
$Z(\mathfrak{g},\mathfrak{k})$ is an associative algebra and a
bimodule over $R(\mathfrak{h})$ \cite{JL}.

Choose linear independent elements $v_1,...,v_n\in \mathfrak{g}$,
 such that  $<v_1,...,v_n>\oplus\mathfrak{ k}=\mathfrak{g}$ as linear spaces over  $\mathbb{C}$. Put $z_i=pv_i mod J'$.
It can be proved that monomials
$\check{z}_1^{m_1}...\check{z}_n^{m_n}$, $m_i\in\mathbb{Z}^+$, form
a bases of $Z(\mathfrak{g},\mathfrak{k})$ over $R(\mathfrak{h})$.

In the case $Z(\mathfrak{o}_{2n+1},\mathfrak{o}_{2n-1})$  put
$$\check{z}_{i\pm n}=pF_{i,\pm n}modJ', \,\,\,\,i=-n,...,n.$$ Notations are taken from \cite{Molev}. There exists an obvious
symmetry $\check{z}_{ij}=\check{z}_{-j-i}$. From previous
considerations  it follows that
$Z(\mathfrak{o}_{2n+1},\mathfrak{o}_{2n-1})$  is generated by
elements  $\check{z}_{ia}$, $i=0,...,n$, $a=\pm n$ or
$\check{z}_{ai}$, $i=0,...,n$, $a=\pm n$.

Sometimes it is more useful to use the generators

$$z_{i\pm n}=\check{z}_{i\pm n}(f_i-f_{i-1})...(f_i-f_{-n+1}),$$

where $$f_i=F_{ii}+ \rho_i, \text{ for $i>0$
},\,\,\,f_0=-\frac{1}{2},\,\,\, f_{-i}=-f_i,$$ and
$$\rho_i=i-\frac{1}{2} \text{ for $i>0$ and } \rho_{-i}=-\rho_{i}.$$

 In particular

 $$z_{0n}=\check{z}_{0n}\prod_{i=1}^{n-1}(F_{ii}+i-\frac{1}{2}).$$

The Mickelsson-Zhelobenko algebra
$Z(\mathfrak{o}_{2n+1},\mathfrak{o}_{2n-1})$  acts on the space
$\oplus_{\mu}V_{\mu}^+$ (see \cite{Molev}). A weight $\mu$ changes
under this action according to the following rule. Let $a$ be $\pm
n$ and
$\mu+\delta_i=(\mu_1,...,\mu_{i-1},\mu_i+1,\mu_{i+1},...,\mu_{n-1})$.
Then for $i=1,...,n-1$ the following holds

$$z_{ia}:V_{\mu}^+\rightarrow V_{\mu+\delta_i}^{+},\,\,\,\,z_{ai}:V_{\mu}^+\rightarrow V_{\mu-\delta_i}^{+}$$

Elements $z_{0a}$ do not change a $\mathfrak{o}_{2n-1}$-weight, that
is they map $V_{\mu}^+$ into itself.

The pfaffians $PfF_{\hat{n}}$, $PfF_{\widehat{-n}}$, as it was
pointed out in the lemma \ref{hw}, also  act on each
 space $V_{\mu}^+$.

Obviously the images  $pPfF_{\widehat{n}}modJ'$,
$pPfF_{\widehat{-n}}modJ'$ of pfaffians in the Mickelsson-Zhelobenko
algebra act in the same way as the corresponding pfaffians. In the
next section it is proved that
$pPfF_{\widehat{n}}modJ'=C\check{z}_{n0}$, where $C\in
U(h_{\mathfrak{o}_{2n-1}})$. The element $C$ is calculated
explicitly.

\subsection{Images of pfaffians in the Mickelson-Zhelobenko alebra I}
\label{mjp1}


\begin{defn}
 A product  of root and Cartan elements in the universal enveloping
  algebra is called normally ordered if in it at first (from the left side)  the negative root elements
   occur, then Cartan elements occur and at the end  positive root elements occur.
\end{defn}

  Every product of  root and Cartan elements  equals to a sum of normally ordered products.


\begin{prop} Let  $I\subset\{-n+1,...,n-1\}$ be a subset which is not symmetric with respect to zero. Then $pPfF_I=0$ in
$Z(\mathfrak{o}_{2n+1},\mathfrak{o}_{2n-1})$ or in
$Z(\mathfrak{o}_{2n},\mathfrak{o}_{2n-2})$.
\end{prop}
 \proof
According to the definition a pfaffian a sum over permutations.  The
summands  are products of root vectors and Cartan elements  of
$\mathfrak{o}_N$

 The sum of root corresponding to element of each product equals
 $-\sum_{i\in I}e_i$. Since the set $I$ is nonsymmetric one has $-\sum_{i\in I}e_i\neq 0$.

Impose a normal ordering in every summand. When one does the normal
ordering new summands appear. But from the equality
$[e_{\alpha},e_{\beta}]=N_{\alpha,\beta}e_{\alpha+\beta}$ it follows
that the sum of roots corresponding to the elements of these new
products is again  $-\sum_{i\in I}e_i$.

Since $-\sum_{i\in I}e_i\neq 0$  in every normally ordered summand in the pfaffian there is a root element.
 These elements  either are zero modulo $J'$,  if there is  a positive root element, or become vanish  after multiplication by $p$,
  if there is a negative root element.

\endproof

Let us give formula for the image of a pfaffian whose indexing set
$I$ is symmetric and is contained in $\{-n+1,...,n-1\}$. In this
case the calculation of the image in the Mickelsson-Zhelobenko
algebra is equivalent to the calcualtion of the image of the
pfaffian under the Harish-Candra homomorphism. This calculation was
done in \cite{ce2} (proposition 7.1), the result is the following.

\begin{prop}\cite{ce2}
\label{tps}
$PfF_I=\frac{1}{(\frac{k}{2})!}D_{\frac{k}{2}}(F_{i_1i_1},...,F_{i_{\frac{k}{2}\frac{k}{2}}}),$

where $D_r(h_1,...,h_r)=\prod_{i=1}^r(h_i-\frac{r}{2}+i)$

\end{prop}

\subsection{Images of pfaffians in the Mickelsson-Zhelobenko algebra II}
\label{mjp2}


In the previous subsection the images in  $Z(\mathfrak{o}_{2n+1},\mathfrak{o}_{2n-1})$
 or
$Z(\mathfrak{o}_{2n},\mathfrak{o}_{2n-2})$ of pfaffians $PfF_I$ were
found, where $I\subset\{-n+1,...,\widehat{0},...,n-1\}$.

Now let us found an image in
$Z(\mathfrak{o}_{2n+1},\mathfrak{o}_{2n-1})$ of the pfaffian
$PfF_{\widehat{n}}$.

To formulate the next theorem define a polynomial $C_n$.

\begin{defn}
\label{ch} Let
$C_{n-1}(h_1,...,h_{n-1})=(-1)^{n-1}D_{n-1}(h_1,...,h_{n-1})-\\-4\sum_{i=1}^{n-1}(-1)^{t+1}D_{n-2}(h_1,...,\widehat{h_i},...,h_{n-1})$

\end{defn}

\begin{thm}
\label{tm} The image of $PfF_{\widehat{n}}$ in
$Z(\mathfrak{o}_{2n+1},\mathfrak{o}_{2n-1})$
 equals $\check{z}_{n0}C_{n-1}(F_{11},...,F_{(n-1)(n-1)})$.

\end{thm}
\proof

Take a set of indices of type
$I=\{-n,-i_{\frac{k-2}{2}},...,-i_{1},0,i_1,...,i_{\frac{k-2}{2}}\}$



By Lemma \ref{minorn} the following equality takes place

$$PfF_{I}=\sum_{i\in I\setminus \{-n\}}\sum_{I'\sqcup
I''=I\setminus\{i,-n\}}\frac{(\frac{|I'|}{2})!(\frac{|I''|}{2})!}{(\frac{k}{2})!}(-1)^{(I'-niI'')}PfF_{I'}F_{ni}PfF_{I''}.$$

To find  the image in the Mickelson-Zhelobenko algebra of the sum
$\sum_{I'\sqcup I''=I\setminus\{i,-n\}}$ divide the summands into
three groups: $1)$ those for which $i=0$, $2)$ those for which
$i<0$, $3)$ those for which $i>0$.

Let us found the  image of summands for which $i=0$. In this case
$PfF_{I'}$ and $PfF_{I'}$ commute with $F_{n0}$. Note that
$(-1)^{(I'-n0I'')}=(-1)^{(I'I'')}(-1)^{\frac{k}{2}-1}$ (to prove
this firstly move $-n,0$  to two first places and then move $0$ to
the right place, then the signs $(-1)^{|I'|}$, $(-1)^{|I'|-1}$,
$(-1)^{\frac{k}{2}}$ appear).

Using the corollary \ref{corl2} one gets that the image sum these
summands equals
\begin{center}
$(-1)^{\frac{k}{2}-1}\sum_{I'\sqcup
I''=I\setminus\{-n,0\}}\frac{(\frac{|I'|}{2})!(\frac{|I''|}{2})!}{(\frac{k}{2})!}(-1)^{(I'I'')}PfF_{I'}PfF_{I''}F_{n0}=
\frac{1}{\frac{k}{2}}(-1)^{\frac{k}{2}-1}PfF_{I\setminus\{-n0\}}F_{n0}=(-1)^{\frac{k}{2}-1}PfF_{I\setminus\{-n,0\}}F_{n0}.$
\end{center}

Since the sets of indices $\pm (I\setminus\{-n0\})$ and $\{-n,0\}$
do not intersect, one can apply the projector $p$ and equivalence
$modJ'$ to each multiple.

Thus  the image of these summands is
$$\check{z}_{n0}(pPf_{I\setminus\{-n0\}}mod J').$$

Found the image of summands
$$\frac{(\frac{|I'|}{2})!(\frac{|I''|}{2})!}{(\frac{k}{2})!}(-1)^{(I'-niI'')}PfF_{I'}F_{ni}PfF_{I''},$$
for which $i\neq 0$. Let  $i>0$. Then change $F_{ni}$ and
$PfF_{I''}$. One obtains an expression
$$\frac{(\frac{|I'|}{2})!(\frac{|I''|}{2})!}{(\frac{k}{2})!}(-1)^{(I'-niI'')}(PfF_{I'}PfF_{I''}F_{ni}-PfF_{I'}[PfF_{I''},F_{ni}]).$$
 Now let $i<0$.
Change $F_{ni}$ and $PfF_{I'}$,
 one gets

 $$\frac{(\frac{|I'|}{2})!(\frac{|I''|}{2})!}{(\frac{k}{2})!}(-1)^{(I'-niI'')}(F_{ni}PfF_{I'}PfF_{I''}-[PfF_{I'},F_{ni}]PfF_{I''}).$$

Consider the case  $i>0$. In the last expression the first summand
has a zero image in the Mickelsson-Zhelobenko algebra by the
following reason.
 The sum of roots corresponding to the elements $F_{ij}$ that
 participate in the expression for
$PfF_{I'}F_{ni}PfF_{I''}$ equals to $e_{n}$. The element $F_{ni}$
corresponds to the root $e_{n}-e_{i}$. Thus the sum of roots
corresponding to the elements  $PfF_{I'}PfF_{I''}$ equals $-e_i$.
 Express $PfF_{I'}PfF_{I''}$ as a sum of normally ordered products. Since $i>0$ than in every obtained
   normally product there is a negative root element of the algebra  $\mathfrak{o}_{2n-1}$. Thus after applying the extremal projector $p$ the expression
$PfF_{I'}PfF_{I''}$ vanishes.

In the case $i<0$ it is similarly proved that the first summand has
a zero image in the Mickelsson-Zhelobenko algebra.

Now consider the second summand $$-PfF_{I'}[PfF_{I''},F_{ni}]$$ in
the case $i>0$ or $$-[PfF_{I'},F_{ni}]PfF_{I''}$$ in the case $i<0$.
In the first case if $-i\notin I''$ it is zero and it equals to
$-PfF_{I'}PfF_{I''}\mid_{-i\mapsto -n}$ otherwise. In the second
case   if $-i\notin I'$ it is zero and it equals  to
$-PfF_{I'}\mid_{-i\mapsto -n}PfF_{I''}$ otherwise.

Thus the image of summands for which $i\neq 0$ equals to the image
of the expression
\begin{center}
$ -\sum_{i\in I\setminus \{-n\},i>0}\sum_{I'\sqcup
I''=I\setminus\{-n,i\},-i\in
I''}\frac{(\frac{|I'|}{2})!(\frac{|I''|}{2})!}{(\frac{k}{2})!}(-1)^{(I'-niI'')}PfF_{I'}PfF_{I''}\mid_{-i\mapsto
-n}- \sum_{i\in I\setminus \{-n\},i<0}\sum_{I'\sqcup
I''=I\setminus\{-n,i\},-i\in
I'}\frac{(\frac{|I'|}{2})!(\frac{|I''|}{2})!}{(\frac{k}{2})!}(-1)^{(I'-niI'')}PfF_{I'}\mid_{-i\mapsto
-n}PfF_{I''}$\end{center}

Let us prove a proposition.

\begin{prop}
The expression above equals
 \begin{equation}\label{gy}-2\sum_{t=-\frac{k}{2},\neq
0}^{\frac{k}{2}}(-1)^{\frac{k}{2}-t-1}\sum_{J'\sqcup
J''=I\setminus\{\pm
i\}}\frac{(\frac{|J'|}{2})!(\frac{|J''|}{2})!}{(\frac{k}{2})!}(-1)^{(J'J'')}PfF_{J'}PfF_{J''}.\end{equation}

\end{prop}

\proof To prove this let us  firstly calculate the sign
$(-1)^{(I'-niI'')}$. The sign $(-1)^{(I'-niI'')}$  differs from the
sign $(-1)^{(I'I'')}$ by the sign of the permutation which moves
$-n,i$ to their right places. This permutation can be done as
follows: first of all move $-n,i$ to two last places,  then move $i$
to it's right place. If $i=i_t$, then
$$(-1)^{(I'-niI'')}=(-1)^{(I'I'')}(-1)^{(|I'|+|I'|+\frac{k-2}{2}-t)}=(-1)^{\frac{k}{2}-t-1}(-1)^{(I'I'')}.$$

Secondly compare $PfF_{I'}\mid_{-i\mapsto -n}$,
$PfF_{I''}\mid_{-i\mapsto -n}$ and
$PfF_{(I'\setminus\{-i\})\cup\{-n\}}$,
$PfF_{(I''\setminus\{-i\})\cup\{-n\}}$ respectively. Here it is
assumed that $i\in I'$ and $i\in I''$.
 In all these expressions at first, the index  $-i$  is changed to
$n$,  but then in the last two expressions the new set of indices is naturally ordered.
  Thus $PfF_{I'}\mid_{-i\mapsto -n}$ and  $PfF_{(I'\setminus\{-i\})\cup\{-n\}}$, $PfF_{I''}\mid_{-i\mapsto -n}$ and  $PfF_{(I''\setminus\{-i\})\cup\{-n\}}$, differ by the sign of this ordering.

  For summands in the sum $\sum_{i\in I\setminus \{-n\},i<0}\sum_{I'\sqcup
I''=I\setminus\{-ni\},-i\in I''}$ denote
  $$J':=(I'\setminus\{-i\})\cup\{-n\}, \,\,\,J'':=I''.$$ One obtains
  $$(-1)^{(I'I'')}PfF_{I'}\mid_{-i\mapsto
  -n}PfF_{I''}=(-1)^{(J'J'')}PfF_{J'}PfF_{J''}.$$ The sign that appears after the ordering is contained in
  $(-1)^{(J'J'')}$.

  Analogously for the summands in the sum $\sum_{i\in I\setminus \{-n\},i<0}\sum_{I'\sqcup
I''=I\setminus\{-ni\},-i\in I''}$, denote
$$J':=I',\,\,\,J'':=(I''\setminus\{-i\})\cup\{-n\}.$$ One obtains that
$$(-1)^{(I'I'')}PfF_{I'}\mid_{-i\mapsto
  -n}PfF_{I''}=(-1)^{(J'J'')}PfF_{J'}PfF_{J''}.$$

In both cases one has $J'\sqcup J''=I\setminus\{\pm i \}$. Also $|J'|=|I'|$, $|I''|=|J''|$.

Note that a pair of sets $J',J''$ occurs twice.  First  as
$(I'\setminus\{-i\})\cup\{-n\}$, $I''$, second as $I'$,
$(I''\setminus\{-i\})\cup\{-n\}$.

Thus one obtains that the considered sum of images of summands for
which $i\neq 0$ is given by the expression \ref{gy}

 $$2\sum_{t=-\frac{k}{2},\neq
0}^{\frac{k}{2}}(-1)^{\frac{k}{2}-t-1}\sum_{J'\sqcup
J''=I\setminus\{\pm
i_t\}}\frac{(\frac{|J'|}{2})!(\frac{|J''|}{2})!}{(\frac{k}{2})!}(-1)^{(J'J'')}PfF_{J'}PfF_{J''}$$

The proposition is proved.
\endproof

This expression \ref{gy} equals
$$2\sum_{t=-\frac{k}{2},\neq
0}^{\frac{k}{2}}(-1)^{\frac{k}{2}-t-1}PfF_{I\setminus\{\pm i_t\}}
=4\sum_{t=1}^{\frac{k}{2}}(-1)^{\frac{k}{2}-t-1}PfF_{I\setminus\{\pm
i_t\}} $$ (Corollary \ref{corl2}).

Finally one has \begin{equation} \label{eee}
PfF_{I}=(-1)^{\frac{k}{2}-1}\check{z}_{n0}(pPfF_{I\setminus
\{-n,0\}}mod
J')-4\sum_{t=1}^{\frac{k}{2}}(-1)^{\frac{k}{2}-t-1}PfF_{I\setminus\{\pm
i_t\}}\end{equation}

Note that $PfF_{I\setminus\{-n,\pm i_t\}}$  is a pfaffian
$PfF_{I^{t}}$  for a new indexing set $I^t=I\setminus\{\pm i_t\}$.

This set is of the same type as $I$.  Apply to   each pfaffian
$PfF_{I^{t}}$   the equality (\ref{eee}). For each $t$ there appears
a summand
$$(-1)^{\frac{k-2}{2}-1}\check{z}_{n0}pPfF_{I^t\setminus
\{-n0\}}=(-1)^{\frac{k}{2}-2}\check{z}_{n0}pPfF_{I\setminus \pm
i,0,-n}.$$

Also there appear summands  $$ PfF_{I^t\setminus\{\pm i_s\}}=\pm
PfF_{I\setminus\{\pm i_t,\pm i_s\}}.$$  But the sum of these
summands over $t$ and $s$ is zero. Let $0<t<s$. If this summand
comes from the summand $ PfF_{I\setminus\{\pm i_s\}}$ in
(\ref{eee}), then it appears with the sign
$(-1)^{\frac{k}{2}-s-1}(-1)^{(\frac{k}{2}-1)-t-1}$.  If it comes
from the summand $ PfF_{I\setminus\{\pm i_t\}}$ in (\ref{eee})  then
it has the sign
$(-1)^{\frac{k}{2}-t-1}(-1)^{(\frac{k}{2}-1)-(s-1)-1}$. The sum of
these signs is zero.

Hence
$$PfF_{I}=(-1)^{\frac{k}{2}-1}\check{z}_{n0}pPfF_{I\setminus
\{-n,0\}}-4\sum_{t=1}^{\frac{k}{2}}(-1)^{\frac{k}{2}-t-1}(-1)^{\frac{k}{2}-2}\check{z}_{n0}pPfF_{I\setminus
\pm i,0,-n}.$$

Apply the obtained formulae  to $I=\hat{n}$. Recall that according
to Proposition \ref{tps} one has $pPfF_{\widehat{ 0, \pm
n}}=D_{n-1}(F_{11},...,F_{(n-1)(n-1)})$, and $pPfF_{\widehat{ \pm
i,0, \pm
n}}=D_{n-2}(F_{11},...,\widehat{F_{ii}},...,F_{(n-1)(n-1)})$. Thus
one proves Theorem.
\endproof




\section{ Action of pfaffians in the multiplicity space and on the Gelfand-Tsetlin-Molev base.}
\label{dpb}

Using Theorem \ref{tm} obtain formulaes of the action of
$PfF_{\widehat{n}}$ on a base in $V_{\mu}^{+}$ and then on the
Gelfand-Tsetlin-Molev bases in $V$.

Everywhere below indices $a,b$  belong to the set $\{-n,n\}$. We
have introduced the notation $\rho_i=i-\frac{1}{2}$ for $i>0$, and
also $\rho_{-i}=-\rho_{i}$. Also we have denoted
$f_i=F_{ii}+\rho_{i}$ for $i>0$, $f_0=\frac{1}{2}$  and
$f_{-i}=-f_{i}$. Introduce a new notation $$g_i=f_{i}+\frac{1}{2}
\text{ for all } i.$$

Define elements $Z_{ab}(u)$ of the Mickelson-Zhelobenko algebra by the formulae (see $\S 9.3$ in \cite{Molev})

$$Z_{ab}(u)=-(\delta_{ab}(u+\rho_n+\frac{1}{2})+F_{ab})\Pi_{i=-n+1}^{n-1}(u+g_i)+\sum_{i=-n+1}^{n-1}z_{ai}z_{ib}\Pi_{j=-n+1,j\neq
i}^{n-1}\frac{u+g_i}{g_i-g_j}.$$


The mapping $s_{ab}(u)\mapsto u^{-2n}Z_{ab}(u)$ defines a homomorphism of the twisted yangian $Y(\mathfrak{o}_2)\rightarrow
Z(\mathfrak{o}_{2n+1},\mathfrak{o}_{2n-1})$ (see $\S 9.3$ in
\cite{Molev}).

Since the Mickelson-Zhelobenko algebra
$Z(\mathfrak{o}_{2n+1},\mathfrak{o}_{2n-1})$ acts on
$\oplus V_{\mu}^{+}$, the space $\oplus_{\mu} V_{\mu}^{+}$ becomes a
$Y(\mathfrak{o}_2)$-representation. It can be easily proved that the defined action of the yangian preserves the
$\mathfrak{o}_{2n-1}$-weights and hence each space $V_{\mu}^{+}$ is a $Y(\mathfrak{o}_2)$-representation.
 This representation is a sum of two irreducible $U,U'$. The types of $U,U'$ are known.

Let the highest weight of $V$ be
$\lambda=(\lambda_1,...,\lambda_n)$, where
$$0\geq\lambda_1\geq...\geq\lambda_n.$$ A base of the
$Y(\mathfrak{o}_2)$-module $V_{\mu}^{+}$ explicitly is constructed
as follows (all facts and notations are taken from $\S 9.6$ in
\cite{Molev}).

  One takes a collection of numbers $(\sigma,\nu_1,...,\nu_n)$,
satisfying the following conditions
\begin{equation}
\label{nu1}
0\geq
\nu_1\geq\lambda_1\geq\nu_2\geq\lambda_2\geq...\geq\lambda_{n-1}\geq\nu_n\geq
\lambda_n
\end{equation}
\begin{equation}
\label{nu2}
0\geq
\nu_1\geq\mu_1\geq\nu_2\geq\mu_2\geq...\geq\mu_{n-1}\geq\nu_n
\end{equation}

The numbers $\nu_i$ are integers if $\lambda_i$ are integers, and
$\nu_i$ are half integers, if $\lambda_i$ are half integers.  The
number $\sigma$  equals $0$ or $1$, if $\nu_1\neq 0$, and $\sigma$
equals $0$ if $\nu_1=0$.

Put $$\gamma_i=\nu_i+\rho_i+\frac{1}{2}.$$

Let $\xi$ be a highest weight vector of the $\mathfrak{o}_{2n+1}$-module $V$.
 There exists a base of the  $Y(\mathfrak{o}_2)$-module $V_{\mu}^{+}$ formed by
 vectors

$$\xi_{\sigma,\nu}=z_{no}^{\sigma}\Pi_{i=1}^{n-1}z_{ni}^{\nu_i-\mu_i}z_{i-n}^{\nu_i-\lambda_i}\Pi_{k=l_n}^{\gamma_n-1}Z_{n-n}(k)\xi,$$
 where $$l_n=\lambda_n+\rho_n+\frac{1}{2}.$$

Put $\bar{\sigma}=\sigma+1\,\,mod 2$

Write an action of $z_{no}$ on these vectors  following $\S 9.6$ in \cite{Molev}.

If $\sigma=0$, then
$$z_{n0}\xi_{\sigma,\nu}=\xi_{\bar{\sigma},\nu}.$$

If $\sigma=1$, then
$$z_{n0}\xi_{\sigma,\nu}=(-1)^n\sum_{j=1}^n\Pi_{t=1,t\neq
j}^n\frac{-\gamma_t^2}{\gamma_j^2-\gamma_t^2}\xi_{\bar{\sigma},\nu+\delta_j}.$$

From here one immediately obtains formulaes of the action of the
pfaffian $PfF_{\widehat{n}}$ on $V_{\mu}^{+}$. These formulaes are
corollaries of Theorem \ref{tm} and the relation between $z_{0n}$
and $\check{z}_{0n}$.

\begin{lem}
If $\sigma=0$, than
$$PfF_{\widehat{n}}\xi_{\sigma,\nu}=C\xi_{\bar{\sigma},\nu}.$$

If $\sigma=1$, than
$$PfF_{\widehat{n}}\xi_{\sigma,\nu}=(-1)^nC\sum_{j=1}^n\Pi_{t=1,t\neq
j}^n\frac{-\gamma_t^2}{\gamma_j^2-\gamma_t^2}\xi_{\bar{\sigma},\nu+\delta_j}.$$

Here
$C=\frac{C_n(\mu_{1},...,\mu_{n-1})}{\prod_{i=1}^{n-1}(\mu_i+i-1)}$
(see the definition \ref{ch}).

\end{lem}

A base in a $\mathfrak{o}_{2n+1}$-module $V$  with the highest
weight $(\lambda_1,...,\lambda_n)$ is constructed inductively by $n$
using the equality $V=\sum_{\mu}V_{\mu}^+\otimes V(\mu)$, where
$V(\mu)$ is a $\mathfrak{o}_{2n-1}$-representation with the highest
weight $\mu$. The result is the following.

Base vectors of $V$ are indexed by tables  $\Lambda$ of type

$\sigma_n,\lambda_{n,1},\lambda_{n,2},...,\lambda_{n,n}$

$\lambda'_{n,1},\lambda'_{n,2},...,\lambda'_{n,n}$

$\sigma_{n-1},\lambda_{n-1,1},\lambda_{n-1,2},...,\lambda_{n-1,n-1}$

...

$\sigma_1,\lambda_{11}$

$\lambda'_{11}$

The restrictions on these numbers are the following:
\begin{enumerate}

\item $\lambda_{ni}=\lambda_i$

\item $\sigma_k=0,1$

\item The equalities hold:

$\lambda'_{k1}\geq\lambda_{k1}\geq\lambda'_{k2}\geq...\geq\lambda'_{k,k-1}\geq\lambda'_{kk}\geq\lambda_{kk}$
when $k=1,...,n$.

$\lambda'_{k1}\geq\lambda_{k-1,1}\geq\lambda'_{k2}\geq...\geq\lambda'_{k,k-1}\geq\lambda_{k-1,k-1}\geq\lambda'_{kk}$
when $k=2,...,n$.

\item If $\lambda'_{k1}=0$, then $\sigma_k=0$.
\end{enumerate}

Derive formulaes for the action of the pfaffian $PfF_{\widehat{n}}$.
Write the equality $V=\sum_{\mu}V_{\mu}^+\otimes V(\mu)$. From one
hand the action of the pfaffian on $V_{\mu}^+$ is already described.
From the other hand the pfaffian commutes with
$\mathfrak{o}_{2n-1}$. Thus the action on
$V=\sum_{\mu}V_{\mu}^+\otimes V(\mu)$
 is written as
$\sum_{\mu}(PfF_{\widehat{n}}\mid_{V_{\mu}^+})\otimes id$. Hence the
pfaffian changes only the two upper rows of the table $\Lambda$
according to the rule  described above.

Write the table $\Lambda$ as $(\sigma,\lambda,\nu,\Lambda')$, where
$\sigma=\sigma_n$, $\nu=\{\lambda_{n}\}$, $\nu=\{\lambda'_{n}\}$ and
$\Lambda'$ is the rest part of the table $\Lambda$. The base vector
corresponding to a table $\Lambda$ we denote as  $\xi_{\Lambda}$ or
$\xi_{\sigma,\nu,\Lambda'}$.

The following theorem is proved

\begin{thm}
\label{maint}

On the vector $\xi_{\Lambda}$ the pfaffian $PfF_{\widehat{n}}$ acts
as follows.

Let $\Lambda=(\sigma,\lambda,\nu,\Lambda')$, where $\sigma,\lambda$ is  first row of  $\Lambda$, $\lambda'$ is the second row of $\Lambda$ and $\Lambda'$ is the rest part of $\Lambda$.

If $\sigma=0$, then
$$PfF_{\widehat{n}}\xi_{\sigma,\nu,\Lambda'}=C\xi_{\bar{\sigma},\nu,\Lambda'}.$$

If $\sigma=1$, then
$$PfF_{\widehat{n}}\xi_{\sigma,\nu,\Lambda'}=(-1)^nC\sum_{j=1}^n\Pi_{t=1,t\neq
j}^n\frac{-\gamma_t^2}{\gamma_j^2-\gamma_t^2}\xi_{\bar{\sigma},\nu+\delta_j,\Lambda'}.$$

Here
$C=\frac{C_n(\lambda_{n-1,1},...,\lambda_{n-1,n-1})}{\prod_{i=1}^{n-1}(\lambda_{n-1,i}+i-1)}$(see
the definition \ref{ch}).

\end{thm}


\section{Appendix. Pfaffians and tensor representations.}
\label{tens}

In the Appendix an action of a pfaffian in an irreducible tensor
representation is investigated.  Base vectors of such a
representation are encoded by  orthogonal Young tableaus
\cite{proc}. In this section at first an action on tensor products
of vectors of the standard representation is calculated
(propositions \ref{stan},\ref{er},\ref{Pfk},\ref{Pft}). Then Theorem
\ref{Pfvt} giving some information about the action on the base
vectors defined by Young tableaus  is proved.  In this theorem an
image of a base vector is expressed as a linear combination of not
necessarily orthogonal tableaus.  So this theorem does not produce a
formulae of an action of a pfaffian in the bases formed by
orthogonal tableaus.

Let $\lambda=(\lambda_1,...,\lambda_n)$ be a highest weight of a
representation, in this section it is suggested to be integer. To
the  highest weight $\lambda$ there corresponds a Young diagram. Let
$V$ be a standard representation of $\mathfrak{o}_N$. Denote by
$e_i$, $i\in\{-n,...,n\}$  unit base vectors of $V$.

 In the space $V^{\otimes m}$ there exists a subspace $V^{[m]}$
 of traceless tensors. A tensor is traceless if for each pair of indices $1\leq p<q\leq n$
it belongs to the kernel of all contractions $V^{\otimes
n}\rightarrow V^{\otimes (n-2)}$, given by the formulaes
$v_1\otimes...\otimes v_n\mapsto (v_p,v_q)v_1\otimes...\otimes
\widehat{v_p}\otimes ...\otimes \widehat{v_q}\otimes ...\otimes
v_n$, where $(v_p,v_q)=\delta_{p,-q}$  is a scalar product
corresponding to the form $G$ (see sec. \ref{or}).

Denote as $\mathbb{S}_{\lambda}$ a representation obtained from
$V^{\otimes (\sum \lambda_i)}$ by applying the Young symmerizer $c_{\lambda}$,
corresponding to the Young diagramm $\lambda$.

\begin{thm} (see  $\S$ 19.5 in \cite{FH})
Let $V$ be a standart representation of $o_{N}$. The representation
$\mathbb{S}_{[\lambda]}:=V^{[\sum \lambda_i]}\cap
\mathbb{S}_{\lambda} V$ is irreducible and has the highest weight
$\lambda=(\lambda_1,...,\lambda_n)$
\end{thm}

To every Young tableau (a Young diagramm filled by numbers) there
corresponds a vector $v_{T}$ in $\mathbb{S}_{[\lambda]}$. To define
it let us enumerate places of the Young diagram by the numbers
$1,...,m$, where $m=\lambda_1+...+\lambda_n$.


\begin{defn}
\label{vt} Let  in the tableau $T$ on the place $i$ stand the number
$t_i$. Take the tensor $e_{t_1}\otimes ... \otimes e_{t_{m}}$ and
apply to it the Young symmertizer $c_{\lambda}$ corresponding to the
diagram.  Take the projection  of $c_{\lambda}(e_{t_1}\otimes ...
\otimes e_{t_{m}})$ to the space of traceless tensors. Denote the
resulting tensor as $v_{T}$.
\end{defn}

The tensors $v_T$ are not linearly independent. To obtain
independent tensors $v_{T}$ one must take $v_{T}$ corresponding only
to the so called orthogonal Young tableaus.

Let us find an action of a pfaffian $PfF_I$, $|I|=k$ on the vectors  $v_T$ of a $\mathfrak{o}_{N}$-representation, given by Young tableaus.

It is done in several steps. At first step the action on the vectors
$e_r$ of standard representation is described. Then the action on
the vectors  $e_{r_2}\otimes e_{r_4}...\otimes e_{r_{t}}$, where
$t<\frac{k}{2}$ is considered. Then the cases $t=\frac{k}{2}$ and
$t>\frac{k}{2}$ are considered. Using these formulaes the action on
$v_{T}$ is described.

\begin{prop}\label{sta5}

On the base  vectors $e_{-2},e_{-1},e_0,e_1,e_2$ of  the standard
representation of $\mathfrak{o}_{5}$ the pfaffians $PfF_I$ where
$|I|=4$ act as zero operators.
\end{prop}
\proof

The proposition is proved  by direct calculation using the
formulaes, where $a\star b=\frac{1}{2}(ab+ba)$

$PfF_{\widehat{-2}}=F_{0-1}\star F_{-21}-F_{-1-1}\star
F_{-20}+F_{-2-1} \star F_{-10}$

$PfF_{\widehat{-1}}=F_{0-2}\star F_{-21}-F_{-1-2}\star
F_{-20}+F_{-2-2} \star F_{-10}$

$PfF_{\hat{0}}=F_{1-2}\star F_{-21}-F_{-1-2}\star
F_{-2-1}+F_{-2-2}\star F_{-1-1}$

$PfF_{\hat{1}}=F_{1-2}\star F_{-20}-F_{0-2}\star
F_{-2-1}+F_{-2-2}\star F_{0-1}$

$PfF_{\hat{2}}=F_{1-2}\star F_{-10}-F_{0-2}\star
F_{-1-1}+F_{-1-2}\star F_{0-1}$

\endproof

Prove an analog of the previous statement in an arbitrary dimension

\begin{prop}
\label{stan} On the base  vectors $e_{-n},...,e_n$ of the standard
representation of $\mathfrak{o}_{N}$ the pfaffians $PfF_I$ for
$|I|>2$ act as zero operators.
\end{prop}

Put $q=4$, $p=k-4$ in Lemma \ref{minors1}. One has

$$PfF_Ie_j=\sum_{I'\sqcup I''=I,
|I'|=k-4,|I''|=4}\frac{(\frac{p}{2})!(\frac{q}{2})!}{(\frac{k}{2})!}(-1)^{(I'I'')}PfF_{I'}PfF_{I''}e_j.$$

If $j\notin I''$, then obviously  $PfF_{I''}e_j=0$. If  $j\in I''$,
then  using Proposition \ref{sta5} one also obtains
$PfF_{I''}e_j=0$.

\endproof

Let us find an action of a pfaffian of the order  $k$ on a tensor
product of  $<\frac{k}{2}$ vectors,  that is on a tensor product
$e_{r_2}\otimes e_{r_4}\otimes ...\otimes e_{r_t}$, where $t<k$.

\begin{prop}
\label{er} $PfF_Ie_{r_2}\otimes e_{r_4}...\otimes e_{r_t}=0$ where $t<k$
\end{prop}

\proof The following formulae takes place $\Delta
PfF_I=\sum_{I'\sqcup I''=I}(-1)^{(I'I'')}PfF_{I'}\otimes PfF_{I''}$
(Lemma \ref{comult}).

By definition one has $PfF_Ie_{r_2}\otimes e_{r_4}\otimes...\otimes
e_{r_k}=(\Delta^kPfF_I)e_{r_2}\otimes e_{r_4}\otimes...\otimes
e_{r_k}$. Since $t<k$, the comultiplication $\Delta^kPfF_I$ contains
only  summands in which on some place the pfaffian stands whose
indexing set $I$ satisfies  $|I|\geq 4$ (Lemma \ref{comult}). From
Proposition \ref{sta5} it follows that every such a summand acts as
a zero operator. \endproof

Find an action of a pfaffian of the order  $k$ on a tensor product of $\frac{k}{2}$ vector,  that is on the tensor product
$e_{r_2}\otimes e_{r_4}\otimes ...\otimes e_{r_k}$.

\begin{prop}
\label{Pfk}If $\{r_2,r_4...,r_k\}$ is not contained in $ I$, then
$PfF_Ie_{r_2}\otimes e_{r_4}\otimes...\otimes e_{r_k}=0$.

 Otherwise
take a permutation $\gamma$ of $I$, such that
$(\gamma(i_1),\gamma(i_2),...,\gamma(i_k))=(r_1,r_2,r_3,...,r_{k-1},r_k)$.
Then \begin{center}$PfF_Ie_{r_2}\otimes...\otimes
e_{r_k}=(-1)^{\gamma}(-1)^{\frac{k(k-1)}{2}}\sum_{\delta\in
Aut(r_1,r_3,...,r_{k-1})}(-1)^{\delta}e_{-\delta(r_1)}\otimes
e_{-\delta(r_3)}\otimes ...\otimes e_{-\delta(r_{k-3})}\otimes
e_{-\delta(r_{k-1})}.$\end{center}

\end{prop}

\proof By definition one has $$PfF_Ie_{r_2}\otimes
e_{r_4}\otimes...\otimes e_{r_k}=(\Delta^kPfF_I)e_{r_2}\otimes
e_{r_4}\otimes...\otimes e_{r_k}.$$ Applying many times the formulae
for comultiplication one obtains
$$\Delta^{\frac{k}{2}}PfF_I=\sum_{I^1\sqcup...\sqcup
I^k}(-1)^{(I^1...I^k)}PfF_{I^1}\otimes...\otimes PfF_{I^k}.$$ Using
Proposition \ref{er} one gets  that, only the summands for which
$|I^j|=2, j=1,...,k$ are nonzero operators.



Hence the summation over divisions  can be written in the following
manner.

\begin{center}$PfF_Ie_{r_2}\otimes e_{r_4}\otimes...\otimes
e_{r_k}=\frac{1}{2^{\frac{k}{2}}}\sum_{\sigma\in
S_k}(-1)^{\sigma}F_{-\sigma(i_1)\sigma(i_2)}\otimes...\otimes
F_{-\sigma(i_{k-1})\sigma(i_k)}(e_{r_2}\otimes...\otimes
e_{r_k})=\frac{1}{2^{\frac{k}{2}}}\sum_{\sigma\in
S_k}(-1)^{\sigma}F_{-\sigma(i_1)\sigma(i_2)}e_{r_2}\otimes...\otimes
F_{-\sigma(i_{k-1})\sigma(i_k)}e_{r_k}.$\end{center}

Consider the expression $F_{-\sigma(i_1)\sigma(i_2)}e_{r_2}$. This
is $e_{-\sigma(i_1)}$ if $\sigma(i_2)=r_2$, this is
$-e_{-\sigma(i_2)}$ if $\sigma(i_1)=r_2$ and zero otherwise. Thus
the summand is nonzero only if the permutation $\sigma$ satisfies
the following condition. In each pair
$(\sigma(i_{2t-1}),\sigma(i_{2t}))$ either $\sigma(i_{2t-1})=r_{2t}$
or   $\sigma(i_{2t})=r_{2t}$.

Show that one can consider only the permutations $\sigma$ such that
$\sigma(i_{2t})=r_{2t}$, that is the permutations of type
$(\sigma(i_1),r_2,\sigma(i_2),r_3...,\sigma(i_{k-1}),r_k)$. But when
only  summands corresponding to such permutations are considered one must multiply the resulting
sum on $2^{\frac{k}{2}}$.

It is enough to prove that the permutations
$\sigma=(\sigma(i_1),\sigma(i_2)=r_2,\sigma(i_3)...,\sigma(r_k))$
and
$\sigma'=(\sigma(i_2)=r_2,\sigma(i_1),\sigma(i_3)...,\sigma(r_k))$
give the same input.

Remind that the input for $\sigma$ is
$$(-1)^{\sigma}F_{-\sigma(i_1)\sigma(i_2)}e_{r_2}\otimes...\otimes
F_{-\sigma(i_{k-1})\sigma(i_k)}e_{r_k}.$$ One has from one hand that
$F_{-\sigma(i_1)\sigma(i_2)}e_{r_2}=e_{-\sigma(i_1)}$ and from the
other hand
$F_{-\sigma'(i_1)\sigma'(i_2)}e_{r_2}=-e_{-\sigma'(i_2)}=-e_{-\sigma(i_1)}$,
 Also one has $(-1)^{\sigma}=-(-1)^{\sigma'}$.
Thus the inputs corresponding to $\sigma$ and $\sigma'$  are the
same.

Hence one can consider the only the permutations  $\sigma$ of type
$(\sigma(i_1),r_2,\sigma(i_2),r_3...,\sigma(i_{k-1}),r_k)$ but
multiplying the resulting sum on $2^{\frac{k}{2}}$.

For the permutation  $\sigma$ of type
$(\sigma(i_1),r_2,\sigma(i_2),r_3...,\sigma(i_{k-1}),r_k)$  using
the definition of $\gamma$ one gets

\begin{center}$(-1)^{\sigma}F_{-\sigma(i_1)\sigma(i_2)}e_{r_2}\otimes...\otimes
F_{-\sigma(i_{k-1})\sigma(i_k)}e_{r_k}=
(-1)^{(\sigma(i_1)r_2,...,\sigma(i_{k-1})r_k)}e_{-\sigma(i_1)}\otimes
e_{-\sigma(i_3)}\otimes...\otimes
e_{-\sigma(i_k)}=(-1)^{\frac{k(k-1)}{2}}(-1)^{\gamma}(-1)^{\delta}e_{-\delta(r_1)}\otimes
...\otimes e_{-\delta(r_{k-1})}.$\end{center} Here $\delta$ is a
permutation of the set $\{r_1,r_3,...,r_{k-3},r_k\}$.

The equality
$(-1)^{\frac{k(k-1)}{2}}(-1)^{\delta}(-1)^{\gamma}=(-1)^{\sigma}$
was used.


 Taking the summation over all permutations $\delta$,  one gets
 $$PfF_Ie_{r_2}\otimes e_{r_4}\otimes...\otimes e_{r_k}=(-1)^{\frac{k(k-1)}{2}}(-1)^{\gamma}\sum_{\delta\in
Aut(r_1,...,r_{k-1})}(-1)^{\delta}e_{-\delta(r_1)}\otimes...\otimes
e_{-\delta(r_{k-1})}. $$ \endproof

Finally from the formula $PfF_Ie_{r_2}\otimes e_{r_4}\otimes
...\otimes e_{r_t}=(\Delta^tPfF_I)e_{r_2}\otimes
e_{r_4}\otimes...\otimes e_{r_t}$, as in the proof of Proposition
\ref{Pfk}, one gets the formulae of the action on an arbitrary
tensor $e_{r_2}\otimes...\otimes e_{r_t}$.

\begin{prop}
\label{Pft} $PfF_Ie_{r_2}\otimes e_{r_4}...\otimes
e_{r_t}=\sum_{\{j_2,j_4...,j_k\}\subset\{2,4...,t\}}Pf^{j_2,j_4...,j_t}F_Ie_{r_2}\otimes
e_{r_4}\otimes...\otimes e_{r_t}$. Here $Pf^{j_2,j_4...,j_k}F_I$
acts on the tensor multiples with numbers $j_2,j_4...,j_k$. It's
action is described by  Proposition \ref{Pfk}

\end{prop}

Find the action of the pfaffian $PfF_I$ on a base vector $v_T$.

By definition the action $PfF_I$ on $v_T$ is constructed as follows.

\begin{enumerate}
\item To the tensor product $e_{t_1}\otimes ... \otimes e_{t_{m}}$ the Young symmetrizer $c_{\lambda}$  is applied.
\item The projection  on the space of traceless tensors is applyed.
\item The pfaffian $PfF_I$ is applied.

\end{enumerate}

Change the order of operations.

Since  the Young symmerizer and  the projection  on the space of
traceless tensors commute with the action of $\mathfrak{o}_N$ (and
hance with $PfF_I$) one can first apply the pfaffian, than the
symmetrizer and finally the projection.


Using the propositions  \ref{Pft}, \ref{Pfk} one gets the following theorem.

\begin{thm}
\label{Pfvt} The image of  $PfF_Iv_T$ can be found as follows
\begin{enumerate}
\item In all possible ways in the tableau $T$ choose  $\frac{k}{2}$ places in the Young tableau, such that
 on them different indices $r_2,r_4...,r_{k}\in I$ stand.
Find a permutation $\gamma$ of $I$, such that
$(\gamma(i_1),...,\gamma(i_k))=(r_1,r_2,r_3,...,r_{k-1},r_k)$.

Replaced indices $r_2,r_4,...,r_{k}$ in all possible ways onto the
indices $-r_1,-r_3,...,-r_{k-1}$.

Take the alternative sum of tableaus $\sum \pm T'$, the tableau in
this sum are obtained from the initial one by placing
$-r_1,-r_3...,-r_{k-1}$ on the chosen places. The places in the
tableau are ordered and the sign is defined by the placing of
indices $-r_1,-r_3...,-r_{k-1}$  on these places. The resulting sum
is multiplied by $(-1)^{\gamma}(-1)^{\frac{k(k-1)}{2}}$.

  If it is not possible to chose in $T$  places  in which stand different indices $r_2,r_4,...,r_{k}\in I$, than $PfF_Iv_T=0$ .

\item For every $T'$ the tensor $v_{T'}$ is constructed.

\item The sum  $\sum_{T'}v_{T'}$ is taken. The result is $PfF_Iv_T$.

\end{enumerate}
\end{thm}

\end{document}